\documentclass{svmult}

\usepackage{makeidx}         
\usepackage{graphicx}        
\usepackage{multicol}        
\usepackage[bottom]{footmisc}
\usepackage{natbib}
\usepackage{url}

\renewcommand{\PackageWarningNoLine}[2]{}
\usepackage{amsmath}
\usepackage{amsfonts}
\usepackage{amssymb}%

\usepackage{booktabs} 
\usepackage{xcolor}

\newcommand{\bE}{\mathbf{E}}
\newcommand{\bJ}{\mathbf{J}}
\newcommand{\ri}{{\rm i}}
\newcommand{\re}{{\rm e}}
\newcommand{\tin}{\text{ in }}
\newcommand{\bn}{\mathbf{n}}
\newcommand{\bze}{\textbf{0}}
\newcommand{\ton}{\text{ on }}
\newcommand{\beq}{\begin{equation}}
\newcommand{\eeq}{\end{equation}}
\newcommand{\HocO}{H_0(\text{curl};\Omega)}

\newcommand{\Hsub}{H_{\text{sub}}}
\newcommand{\bv}{{\bf v}}
\newcommand{\curl}{\nabla\times}
\newcommand{\bvb}{\overline{{\bf v}}}
\newcommand{\cV}{{\mathcal V}}
\newcommand{\cI}{{\mathcal I}}
\newcommand{\beqs}{\begin{equation*}}
\newcommand{\eeqs}{\end{equation*}}
\newcommand{\bw}{{\bf w}}
\newcommand{\br}{\mathbf{r}}
\newcommand{\bfr}{\mathbf{r}}
\newcommand{\bF}{\mathbf{F}}
\newcommand{\bU}{\mathbf{U}}
\newcommand{\bV}{\mathbf{V}}
\newcommand{\bW}{\mathbf{W}}
\newcommand{\tOmega}{\widetilde{\Omega}}
\newcommand{\R}{\mathbb{R}}
\newcommand{\cT}{{\mathcal T}}

\newcommand{\matrixA}{A}

\newcommand{\cO}{{\mathcal O}}
\newcommand{\Hprec}{H_{\mathrm{prec}}}

\newcommand{\cC}{{\mathcal C}}

\usepackage{soul}

\usepackage[final]{changes}

\definechangesauthor[name=Marcella,color=blue]{MB}
\definechangesauthor[name=report1,color=red]{R1}
\definechangesauthor[name=report1,color=blue]{R2}

\begin{document}

\title*{A two-level {domain-decomposition} preconditioner for the {time-harmonic} Maxwell's equations}
\titlerunning{A two-level DD preconditioner for the {time-harmonic} Maxwell's equations}
\author{
  Marcella Bonazzoli\inst{1}\and
  Victorita Dolean\inst{1,2}\and
  Ivan G. Graham \inst{3} \and 
  Euan A. Spence \inst{3} \and
  Pierre-Henri Tournier \inst{4}}
  
\authorrunning{M. Bonazzoli, V. Dolean, I.G. Graham, E.A. Spence, P.-H. Tournier}  

\institute{
  \inst{1}
  Universit\'{e} C\^{o}te d'Azur, CNRS, LJAD, France,
  \email{marcella.bonazzoli@unice.fr}\\ 
  \inst{2}
  University of Strathclyde, Glasgow, UK,
  \email{Victorita.Dolean@strath.ac.uk}\\ 
  \inst{3} 
  University of Bath, UK,
  \email{I.G.Graham@bath.ac.uk, E.A.Spence@bath.ac.uk}\\
  \inst{4}
  UPMC Univ Paris 06, LJLL, Paris, France,
  \email{tournier@ljll.upmc.fr}
}
%
%
\maketitle

\abstract*{The construction of fast iterative solvers for the indefinite time-harmonic Maxwell's system at mid- to high-frequency is a problem of great current interest. Some of the difficulties that arise are similar to those encountered in the case of the mid- to high-frequency Helmholtz equation. Here we investigate how two-level domain-decomposition preconditioners recently proposed for the Helmholtz equation work in the Maxwell case, both from the theoretical and numerical points of view.} 

\section{Introduction}

{The construction of fast iterative solvers for the indefinite time-harmonic Maxwell's system at \added[id=R2]{mid- to} high-frequency is a problem of great current interest. Some of the difficulties that arise are similar to those encountered in the case of the \added[id=R2]{mid- to} high-frequency Helmholtz equation. Here we investigate how domain-decomposition (DD) solvers recently proposed for the \deleted[id=R2]{high-frequency} Helmholtz equation work in the Maxwell case.}

{The idea of preconditioning discretisations of the Helmholtz equation with discretisations of the corresponding problem with absorption was introduced in \cite{ErVuOo:04}. In \cite{Graham:2016:RRD}, a two-level domain-decomposition method was proposed that uses absorption, along with a wavenumber dependent coarse space correction. Note that, in this method, the choice of absorption is motivated by the analysis in both \cite{Graham:2016:RRD} and the earlier work \cite{Gander:2015:AGH}.}

{Our aim is to extend these ideas to the time-harmonic Maxwell's equations, both from the theoretical and numerical points of view. These results will appear in full in the forthcoming paper \cite{BoDoGrSpTo:17}.}

{Our theory will apply to} the boundary value problem (BVP)
\begin{equation}
\left\{  
\begin{array}{rl}
\nabla\times(\nabla\times\bE) - (k^2 + \ri \kappa) \bE &= \bJ \quad\tin \Omega\\
\bE \times \bn &= \bze \quad\ton \Gamma:= \partial \Omega
\end{array}\quad \quad \right.\label{eq:BVP}
\eeq
where $\Omega$ is a bounded Lipschitz {polyhedron} in $\mathbb{R}^3$ with boundary $\Gamma$
and outward-pointing unit normal vector $\bn$, $k$ is the wave number, and $\bJ$ is the source term. 
{The PDE in \eqref{eq:BVP} is obtained from Maxwell's equations by assuming} that the electric field $\boldsymbol{\mathcal{E}}$ is of the form $\boldsymbol{\mathcal{E}} (\mathbf{x},t) = \Re (\mathbf{E}(\mathbf{x}) \re^{-\ri\omega t})$, where $\omega>0$ is the angular frequency.
The boundary condition in~\eqref{eq:BVP} is called Perfect Electric Conductor (PEC) boundary condition.
The parameter $\kappa$ dictates the absorption/damping in the problem; in the case of a conductive medium, $\kappa=k\sigma Z$, where $\sigma$ is the electrical conductivity of the medium and $Z$ the impedance.
If $\sigma=0$, the solution is not unique {for all $k>0$} but a sufficient condition for existence of a solution is $\nabla\cdot \bJ = 0$.

We will {also give numerical experiments for the BVP \eqref{eq:BVP} where the PEC boundary condition is replaced by 
an impedance boundary condition, i.e.~the BVP}
\begin{equation}
\left\{
\begin{array}{ll}
\nabla\times(\nabla\times\bE) - (k^2 + \ri \kappa) \bE = \bJ &\tin \Omega\\
(\nabla\times \mathbf{E}) \times \mathbf{n}  - \ri \, k  \; \mathbf{n} \times(\mathbf{E}  \times \mathbf{n}) = \mathbf{0} &\ton \Gamma:= \partial \Omega
\end{array}  
\right.\label{eq:BVPimp} 
\eeq
{In contrast to the PEC problem, the solution of the impedance problem is unique for every $k>0$.}
{There is large interest in solving \eqref{eq:BVP} and \eqref{eq:BVPimp} \emph{both} when $\kappa=0$ \emph{and} when $\kappa\neq 0$. We will consider both these cases, in each case constructing preconditioners by using larger values of $\kappa$.} 
\added[id=R1]{Indeed, a higher level of absorption makes the problems involved in the preconditioner definition more ``elliptic'' (in a sense more precisely explained in \cite{BoDoGrSpTo:17}), thus easier to solve. Note that the absorption cannot increase too much, otherwise the problem in the preconditioner is ``too far away'' from the initial problem.}


\section{Variational formulation and discretisation}


Let $\HocO := \{ \bv \in L^2(\Omega), \curl \bv \in L^2(\Omega), \bv\times \bn=\bze\}$. We introduce the $k$-weighted inner product on $\HocO$:  
$$
(\bv,\bw)_{\text{curl},k} \ = \ (\curl \bv, \curl \bw)_{L^2(\Omega)} + k^2 (\bv,\bw)_{L^2(\Omega)}.  
$$
The standard variational formulation of \eqref{eq:BVP} is:
Given $\bJ \in L^2(\Omega)$, $\kappa\in \mathbb{R}$ and $k>0$, find $\bE \in \HocO$ such that 
\beq\label{eq:vf_intro}
a_\kappa(\bE,\bv) = F(\bv) \,\,\text{ for all }\,  \bv \in \HocO,
\eeq
where 
\beq\label{eq:Helmholtzvf_intro}
a_\kappa(\bE,\bv) := \int_\Omega \curl  \bE\cdot \overline{\curl \bv}  - (k^2 + \ri \kappa) \int_\Omega \bE\cdot\bvb
\eeq
and 
$F(\bv) := \int_\Omega \bJ\cdot\bvb.$
{When} $\kappa > 0$, {it is well-known that the sesquilinear form is coercive (see, e.g., \cite{BoDoGrSpTo:17} and the references therein) and so existence and uniqueness follow from the Lax--Milgram theorem.}



N\'ed\'elec edge elements 
are particularly suited for the approximation of electromagnetic fields. 
They provide a conformal discretisation of $H(\text{curl},\Omega)$, since their tangential component across faces shared by adjacent tetrahedra of a simplicial mesh $\mathcal{T}^h$ is continuous. 
We therefore 
define our approximation space $\cV^h\subset \HocO$ as the {lowest-order} edge finite element space on the mesh $\mathcal{T}^h$ with functions whose tangential trace is zero on $\Gamma$.
More precisely, over each tetrahedron $\tau$, we write the discretised field as $\mathbf{E}_h= \sum_{e \in \tau} c_e {\bf w}_e$, a linear combination with coefficients $c_e$ of the basis functions ${\bf w}_e$ associated with the edges $e$ of $\tau$, and the coefficients $c_e$ will be the unknowns of the resulting linear system. The Galerkin method applied to the variational problem \eqref{eq:vf_intro} is
\beq\label{eq:Galerkin}
\text{find}\,\, \bE_h \in \cV^h\,\, \text{ such that } \,\,a_\kappa(\bE_h,\bv_h) = F(\bv_h) \,\,\text{ for all }\,  \bv_h \in \cV^h.
\eeq
The Galerkin matrix $A_\kappa$ is defined by
$(A_\kappa)_{i j} := a_\kappa(\bw_{e_i}, \bw_{e_j})$
and the Galerkin method is then equivalent to solving the linear system $A_\kappa \bU =\bF$, where $F_i := F(\bw_{e_i})$ and $U_j := c_{e_j}$.


\section{Domain decomposition}


To define appropriate subspaces of  $\cV^h$,  
we start with a collection of open subsets    $\{ \tOmega_\ell: \ell =
1, \ldots, N\}$ of $\R^d$ of maximum diameter $H_{\text{sub}}$  that form an overlapping cover of $\overline{\Omega}$, 
and we set $\Omega_\ell =
\tOmega_\ell \cap \overline{\Omega}$. 
Each $\overline{\Omega}_\ell$  is assumed to be non-empty 
and is assumed to  consist of a union of
elements of the mesh $\cT_h$. Then, for each $\ell = 1, \ldots, N$,  we set 
\beqs
\cV_\ell := \cV^h \cap H_0(\text{curl}, \Omega_\ell),
\eeqs
where $H_0(\text{curl},\Omega_\ell)$ is considered as a subset of $\HocO$ by extending functions in $H_0(\text{curl},\Omega_\ell)$ by zero, thus the tangential 
traces of elements of $\cV_\ell$ vanish on the
internal boundary $\partial \Omega_\ell \backslash \Gamma$ (as well as on $\partial \Omega_\ell \cap \Gamma$). Thus a solve of the Maxwell problem \eqref{eq:vf_intro}  in the space $\cV_\ell$ involves a
PEC boundary condition on $\partial \Omega_\ell$ {(including {any} external parts of $\partial \Omega_\ell$).}
When $\kappa \not = 0$, such solves are always  well-defined {by uniqueness of the solution of the BVP \eqref{eq:BVP}.}

Let $\cI^h$ be the set of interior edges of elements of the triangulation; this set can be identified with the degrees of freedom of $\cV^h$.
Similarly, let $\cI^h(\Omega_{\ell})$ be the set of edges of elements contained in {(the interior of)} $\Omega_{\ell}$ (corresponding to 
degrees of freedom on those edges). We then have that
$\cI^h= \cup_{\ell=1}^N \cI^h(\Omega_\ell)$. For      
$e\in \cI^h({\Omega_\ell})$ and $e' \in \cI^h$, we define the
restriction matrices  
$(R_\ell)_{e, e'} := \delta_{e, e'}$.
{We will assume that we have} matrices $(D_\ell)_{\ell=1}^N$ satisfying 
\begin{equation}
\label{eq:algebPartUnity}
\sum_{\ell = 1}^{N} R_\ell^T D_\ell R_\ell = I;
\end{equation}


\noindent{such} matrices $(D_\ell)_{\ell=1}^N$ are called a \emph{partition of unity}.

For two-level methods we need to define a coarse space.  
Let  $\{\cT^{H}\}$ be  a sequence of shape-regular, tetrahedral meshes on 
$\overline{\Omega}$, with mesh diameter $H$. We assume that each element of $\cT^H$ consists of the union of a set of fine grid elements. 
Let $\cI^H$ be an index set for the coarse mesh edges.
The coarse basis functions $\{{\bf w}^{H}_e\}$ are taken to be N\'ed\'elec edge elements on $\cT^H$ with zero tangential traces on {$\Gamma$.}
From these functions we define 
the coarse space
$
\cV_0 :=  \mathrm{span} \{{\bf w}^H_{e_p}  : p \in \cI^H
\}, 
$
and we define the ``restriction matrix''
\begin{equation}\label{eq:restriction}
(R_0)_{pj} := \psi_{e_{{j}}}(\bw_{e_{{p}}}^{{H}}) {= \int_{e_j} \bw_{e_p}^H \cdot \mathbf{t}},  \quad j \in \cI^h , \quad
p \in \cI^{H},
\end{equation}
where $\psi_{e}$ are the degrees of freedom on the fine mesh.

With the restriction matrices $(R_\ell)_{\ell=0}^{N}$ defined above, 
we define
 $$A_{\kappa,\ell} \ := \ R_\ell  A_\kappa R_\ell ^T, \quad {\ell = 0,\ldots, N}$$
 For $\ell=1,\ldots, N$, the matrix $A_{\kappa,\ell}$ 
is then  just the minor  of $A_\kappa$ corresponding to rows
and columns taken from $\cI^h({\Omega_\ell})$. That is $A_{\kappa,\ell}$ corresponds to the Maxwell  problem on $\Omega_\ell$ with homogeneous PEC boundary condition on $\partial \Omega_\ell \backslash \Gamma$.   The matrix $A_{\kappa,0}$ is 
 the Galerkin matrix for the problem {\eqref{eq:BVP}} discretised in $\cV_0$.
In a similar way as for the global problem it can be proven that matrices $\matrixA_{\kappa,\ell}$, $\ell=0,\ldots, N$, are invertible for all mesh 
sizes $h$ and all choices of  $\kappa \not = 0$. 

In this paper we consider two-level preconditioners, i.e.~those involving both local and coarse solves, except if `1-level' is specified in the numerical experiments.
The classical
\emph{two-level} \emph{Additive Schwarz} {(AS)} and \emph{Restricted Additive Schwarz} {(RAS)} preconditioners for $A_\kappa$ are defined by
\begin{equation}\label{eq:defAS} 
M_{\kappa,\text{AS}}^{-1}: =  \sum_{\ell=0}^N R_\ell^T A_{\kappa,\ell}^{-1} R_\ell\quad \, M_{\kappa,\text{RAS}}^{-1}: =  \sum_{\ell=0}^N R_\ell^T D_\ell A_{\kappa,\ell}^{-1} R_\ell.
\end{equation}
In the numerical experiments we will also consider two other preconditioners: {(i)} $M_{\kappa,\text{ImpRAS}}^{-1}$, which is similar to $M_{\kappa,\text{RAS}}^{-1}$, but the solves with $A_{\kappa,\ell}$ are replaced by solves with matrices corresponding to the Maxwell problem on $\Omega_\ell$ with homogeneous impedance boundary condition on $\partial \Omega_\ell \backslash \Gamma$, 
 and {(ii)} the hybrid version of RAS
\begin{equation}\label{eq:HybRAS} 
M^{-1}_{\kappa,\text{HRAS}}:= (I - \Xi A_{\kappa}) \Biggl(\sum_{\ell=1}^N R_\ell^T D_\ell A_{\kappa,\ell}^{-1} R_\ell \Biggr) (I - A_{\kappa} \Xi) + \Xi,\, \Xi = R_0^T A_{\kappa,0}^{-1} R_0.
\end{equation}
In a similar manner we can define {$M^{-1}_{\kappa,\text{HAS}}$, $M^{-1}_{\kappa,\text{ImpHRAS}}$,} the hybrid versions of AS and ImpRAS.


\section{{Theoretical results}}


{The following result is the Maxwell-analogue of the Helmholtz-result in \cite[Theorem 5.6]{graham:2015:domain} and {appears in \cite{BoDoGrSpTo:17}}.}
We state a version of this result for $\kappa\sim k^2$, but note that \cite{BoDoGrSpTo:17} contains a more general result that, in particular, allows for smaller values of the absorption $\kappa$.

\begin{theorem}[GMRES convergence for left preconditioning {with $\kappa\sim k^2$}] \label{cor:final2}
{Assume that $\Omega$ is a convex polyhedron.}
Let $C_k$ be the matrix representing  the $(\cdot, \cdot)_{\text{curl},k}$ 
inner product on the finite element space $\cV^h$ in the sense that if $v_h, w_h \in \cV^h$ with coefficient vectors $\bV, \bW$ then 
\begin{equation}\label{eq:ip}
(v_h, w_h )_{\text{curl},k}\  =\ \langle \bV, \bW\rangle_{C_k} . 
\end{equation} 
Consider the weighted GMRES method where the residual is minimised in
the norm induced by $C_k$. Let $\br^m$ denote the $m$th {residual} of GMRES applied to the system $A_\kappa$, left  preconditioned with 
$M_{\kappa,\text{AS}}^{-1}$. Then 
\begin{equation}\label{eq:conv_est}
\frac{\Vert \bfr^m \Vert_{C_k}} { \Vert \bfr^0 \Vert_{C_k} } \
\lesssim   \ \Biggl(1- \biggl(1 + \biggl(\frac{H}{\delta}\biggr)^2\biggr)^{-2}\Biggr)^{m/2} \ ,  
\end{equation}
provided the following condition holds: 
\beq\label{eq:E20rpt}
\max\left\{ k\Hsub,\  kH \right\}  \ \leq \  
\cC_1 \biggl(1 + \biggl(\frac{H}{\delta}\biggr)^2\biggr)^{-1}. 
\eeq
where $\Hsub$ and $H$ are the typical diameters of a subdomain and of the coarse grid, $\delta$ denotes the size of the overlap, {and $\cC_1$ is a constant independent of all parameters.}
\end{theorem}
As a particular example we see that,
provided $\kappa\sim k^2$, $H \sim \Hsub \sim k^{-1}$ and $\delta \sim H$ (``generous overlap"),  then GMRES will
converge with a number of iterations independent of all parameters. 
This property is illustrated in the numerical experiments in the next
section.
{A result analogous to Theorem \ref{cor:final2} 
for right-preconditioning 
appears in \cite{BoDoGrSpTo:17}.}


\section{Numerical results}


In this section we will perform several numerical experiments in a cube domain with PEC boundary conditions (Experiments 1-2) or impedance boundary conditions (Experiments 3-4). 
The right-hand side is given by $\mathbf{J}=[f,f,f]$, where $f = -\exp(-400((x-0.5)^2+(y-0.5)^2+(z-0.5)^2))$.

We solve the linear system with GMRES with right preconditioning, starting with a \emph{random initial guess}, which ensures, unlike a zero initial guess, that all frequencies are present in the error; the stopping criterion, with a tolerance of $10^{-6}$, is based on the relative residual. The maximum number of iterations allowed is $200$. We consider a regular decomposition into subdomains (cubes), the overlap for each subdomain is of size $\mathcal{O}(2h)$ (except in Experiment 1, where we take generous overlap) in all directions. All the computations are done in FreeFem++, an open source domain specific language (DSL) specialised for solving BVPs with variational methods (\url{http://www.freefem.org/ff++/}). The code is parallelised and run on the TGCC Curie supercomputer and the CINES Occigen supercomputer.
We assign each subdomain to one processor.
Thus in our experiments the number of processors increases if the number of subdomains increases.
To apply the preconditioner, the local problems in each subdomain and the coarse space problem are solved with a direct solver (MUMPS on one processor). In all the experiments the fine mesh diameter is $h \sim k^{-3/2}$, which is believed to remove the pollution effect.

In our experiments we will often choose $\Hsub \sim H$ and our preconditioners are thus determined by choices of  $H$ and $\kappa$, which we denote by 
$\Hprec$  
and $\kappa_\text{prec}.$ 
The absorption parameter of the problem to be solved is denoted $\kappa_\text{prob}$.
The coarse grid problem is of size $\sim \Hprec^{-2}$ and there are $\sim \Hprec^{-2}$  
local problems  of size $(\Hprec/h)^2$ (case $\Hsub \sim H$). In the tables of results,  $n$ denotes the size of the system  being solved, $n_{\text{CS}}$ the size of the coarse space, the figures in the tables denote the GMRES iterations corresponding to a given method (e.g.~$\#$AS is the number of iterations for the AS preconditioner), whereas Time denotes the total time (in seconds) \added[id=R2]{including both setup and GMRES solve times}. For some of the experiments we compute (by linear least squares)  the  approximate value of $\gamma$ so that the entries of this column grow with $k^\gamma$. We  also compute  $\xi$ so that the entries of the column grow with $n^\xi$ (here $\xi = \gamma\cdot 2/9$, because  $n \sim  (h^{3/2})^3 = k^{9/2}$).  \\

 
{\bf Experiment 1}. 
The purpose of this experiment is to test  the theoretical result which says that even with AS (i.e. when solving PEC local problems), provided $H \sim \Hsub \sim k^{-1}$, $\delta \sim H$ (generous overlap), $\kappa_\text{prob} = \kappa_\text{prec} = k^2$, the number of GMRES iterations should be bounded as $k$ increases. In Table~\ref{tab:1} we compare three two-level preconditioners: additive Schwarz, restricted additive Schwarz, and the hybrid version of restricted additive Schwarz.
Note that in theory we would expect AS to be eventually robust, although its inferiority compared to the other methods is to be expected 
\cite{Graham:2016:RRD}. 
\begin{table}[t] 
\begin{center} 
\begin{tabular}{|c|c|c|c||c|c|c|}
\hline 
$k$ &  $n$ & $N_\text{sub}$&  $n_{\text{CS}}$ &$\#$AS   & $\#$RAS & $\#$HRAS  \\
\hline 
10 & 4.6 $\times 10^5$ & 1000  & 7.9$\times 10^3$ & 53 & 26& 12 \\
15 & 1.5 $\times 10^6$ & 3375 &2.6$\times 10^4$& 59 & 28& 12\\
20 &1.2 $\times 10^7$& 8000 &6.0$\times 10^4$& 76 & 29& 17 \\
\hline
\end{tabular} 
\caption{$\delta \sim H$ (generous overlap), $H \sim \Hsub \sim k^{-1}$, $\kappa_\text{prob} = \kappa_\text{prec} = k^2$.}    
\label{tab:1}
\end{center} 
\end{table} 

\bigskip
{\bf Experiment 2}. In this experiment (Table~\ref{tab:2}) we set $\kappa_\text{prob} = \kappa_\text{prec} = k^2$ and $H \sim \Hsub \sim k^{-0.8}$ and the overlap is  $\cO(2h)$ in all directions. 
As we are not in the case $\Hprec \sim k^{-1}$ and we do not have generous overlap, we do not expect a bounded number of iterations here. Nevertheless, the method still performs well. Not surprisingly, the best method is ImpHRAS, as better transmission conditions at the interfaces between subdomains are used in the preconditioner. It is important to note that the time is growing very much slower than the dimension of the  problem being solved.  
\begin{table}[t] 
\begin{center} 
\begin{tabular}{|c|c|c|c||c|c|c|}
\hline 
$k $ & $n$ & $N_\text{sub}$&  $n_{\text{CS}}$& $\#$RAS ($\#$HRAS) &   $\#$ImpRAS ($\#$ImpHRAS) & Time ImpHRAS \\
\hline 
10   & $3.4\times 10^5$ & 216 & 1.9$\times 10^3$ & 34 (23) & 27 (20) & 11.0 \\
20   & $7.1\times 10^6$ & 1000 & 7.9$\times 10^3$  &43 (31) & 35 (28) & 42.6  \\
30   & $4.1\times 10^7$ & 3375 &2.6$\times 10^4$  & 47 (34) & 39 (32) &100.9 \\
40   & $1.3\times 10^8$ & 6859 & 5.1$\times 10^4$ & 49 (36) & 42 (35) & 264.5 \\
\hline 
$\gamma$ & 4.5 & & &  &  & 2.23\\
\hline
\end{tabular} 
\caption{$\delta \sim 2h$, $H \sim  H_\text{sub} \sim k^{-0.8} $, $\kappa_\text{prob} = \kappa_\text{prec} = k^2$.} 
\label{tab:2}
\end{center} 
\end{table} 

\bigskip
{\bf Experiment 3} In this case we take $\kappa_\text{prob} = k$. 
Moreover, we take impedance boundary conditions on $\partial \Omega$.  We take $ H \sim  H_\text{sub} \sim k^{-\alpha}$, $\kappa_\text{prec} = k^\beta$, and we use ImpHRAS as a preconditioner.

\begin{table}[t] 
\begin{center} 
\begin{tabular}{cc}
\begin{tabular}{|c||c|c|c|c||c|c|c|c|}
\hline 
& \multicolumn{4}{|c||}{$\alpha = 0.6$} & \multicolumn{4}{c|}{$\alpha = 0.8$} \\
\hline
$k$ & $n$ & $N_\text{sub}$& $n_{\text{CS}}$ &     $\#$2-level  & $n$ & $N_\text{sub}$& $n_{\text{CS}}$&  $\#$2-level   \\
\hline 
10 & $2.6\times 10^5$ & 27 & 2.8$\times 10^2$& 31 & $3.4\times 10^5$ & 216 & $1.8\times 10^3$& 29 \\
20 & $6.3\times 10^6$ & 216 & $1.9\times 10^3$& 87 & $7.1\times 10^6$ & 1000 & $7.9\times 10^3$& 60\\
30 & $3.3\times 10^7$ & 343& $2.9\times 10^3$ & {148} & $4.1\times 10^7$ & 3375 & $2.5\times 10^4$ & 90 \\
40 &$1.1\times 10^8$ & 729& $5.9\times 10^3$ & 200 & $1.3\times 10^8$ & 6859& $5.1\times 10^4$ & 154\\
\hline
\end{tabular} 
\\
\\
\begin{tabular}{|c|c|c|c||c|c|}
\hline 
\multicolumn{4}{|c||}{} & $\beta = 1$ & $\beta = 2$ \\
\hline
$k$ & $n$ & $N_\text{sub}$& $n_{\text{CS}}$ & $\#$2-level(Time) &    $\#$2-level(Time)  \\
\hline 
10 & $3.4\times 10^5$ & 216 & $1.8\times 10^3$ & 29 (12.9) & 37 (13.1)\\
20 & $7.1\times 10^6$ & 1000 & $7.9\times 10^3$ & 60 (63.7) & 70 (69.8) \\ 
30 & $4.1\times 10^7$ & 3375 & $2.5\times 10^4$ & 90 (200.4) & 101 (221.2) \\
40 & $1.3\times 10^8$ & 6859& $5.1\times 10^4$ & 154 (771.7) & 137 (707.6) \\
\hline
$\gamma$ & 4.5  & & 2.4 &1.2 (2.9) &0.94 (2.8)\\
\hline 
$\xi$ & 1.0  & &0.5  &0.3 (0.6) &0.2 (0.6)\\
\hline 
\end{tabular} 
\end{tabular} 
\caption{$\kappa_\text{prob} =k$, $\delta \sim 2h$, $H \sim  H_\text{sub} \sim k^{-\alpha} $, $\kappa_\text{prec} = k^{\beta}$; \replaced[id=R2]{Top:}{Left:} $\beta=2$, $\alpha=0.6,0.8$; \replaced[id=R2]{Bottom:}{Right:} $\alpha=0.8$, $\beta=1,2$.} 
\label{tab:3}
\end{center} 
\end{table} 

In Table~\ref{tab:3} on the {bottom} we see that the dimension of the coarse space is 
$$
n_{\text{CS}} = (k^{-0.8})^{-3} = k^{2.4} = \mathcal{O}(n^{0.5}).
$$  
This is reflected   in the $\gamma$ and $\xi$ figures in the $n_{\text{CS}}$ column.
For this method the reduction factor $n_{\text{CS}}/n$ is substantial (about  
{$3.9 \times 10^{-4} $} when $k = 40$). 
The computation time grows only slightly faster than the dimension of the coarse space, showing (a) weak scaling and (b) MUMPS is still performing close to optimally for Maxwell systems of size $5 \times 10^4$. Iteration numbers are growing with about $n^{0.3}$ at worst. Note that the iteration numbers may be improved by separating the coarse grid size from the subdomain size, making the coarse grid finer and the subdomains bigger.

\bigskip
{\bf Experiment 4}. Here we solve the pure Maxwell problem without absorption, i.e.~$\kappa_\text{prob} = 0$, with impedance boundary conditions 
on $\partial \Omega$. In the preconditioner we take $\kappa_\text{prec} = k$.    
Results are given in Table \ref{tab:4}, {where $H_\text{sub}  \sim k^{-\alpha}$, $H \sim k^{-\alpha'}$}.  These methods are close to being load balanced in the sense that the coarse grid and subdomain problem size are  very similar  when $\alpha + \alpha' = 3/2$.  

Out of the methods tested, the 2-level method (ImpHRAS) with $(\alpha, \alpha')  = (0.6, 0.9)$ gives the best iteration count, but is more expensive. 
The method $(\alpha, \alpha')  = (0.7, 0.8)$  is 
faster 
but its   
iteration count grows more quickly, so its   advantage will  
diminish as $k$ increases further. 
\deleted[id=MB]{We have no explanation for the curious reduction in iterations   in the 2-level method as $k$ increases for $(\alpha, \alpha')  = (0.6, 0.9)$.} 
For $(\alpha, \alpha') = (0.6,0.9)$ the coarse grid size grows with 
$\mathcal{O}(n^{{0.64}}  )$ while the time grows with \replaced[id=MB]{$\mathcal{O}(n^{{0.65}}  )$}{$\mathcal{O}(n^{{0.80}}  )$}. For    
$(\alpha, \alpha') = (0.7,0.8)$ the rates are  $\mathcal{O}(n^{{0.54}}  )$  
and  \replaced[id=MB]{$\mathcal{O}(n^{{0.69}}  )$}{$\mathcal{O}(n^{{0.75}}  )$}.  The subdomain problems are solved on individual processors 
so the number of processors  used grows as $k$ increases. 
In the current implementation a sequential direct solver on one processor is used to factorize the coarse problem matrix, which is clearly a limiting factor for the scalability of the algorithm. The timings could be significantly improved by using a distributed direct solver, or by adding a further level of domain decomposition for the coarse problem solve.

\begin{table}
\begin{center}
\begin{tabular}{|c|c|c||c|c|c|c|c|}
\hline
\multicolumn{3}{|c||}{} & \multicolumn{5}{|c|}{$\alpha=0.6$, $\alpha'=0.9$}\tabularnewline
\hline 
$k$ & $n$ & $N_\text{sub}$ & $\#$2-level & $n_\text{CS}$ & Time & $\#$1-level & Time\tabularnewline
\hline 
10 & $2.6 \times 10^5$ & 27 & 20 & $2.9 \times 10^3$ & 16.2& 37 & 13.7\tabularnewline
15 & $1.5 \times 10^6$ & 125 & 26 & $1.0\times 10^4$ & 25.5& 70 &  26.1\tabularnewline
20 & $5.2 \times 10^6$ & 216 & 29 & $2.1\times 10^4$ & 52.0 & 94 & 60.6\tabularnewline
25 & $1.4\times 10^7$ & 216 & 33 & $4.4\times 10^4$ & 145.5 & 105 &  191.2\tabularnewline 
30 & $3.3\times 10^7$ & 343 & 38 & $6.9\times 10^4$ & 380.4 & 132 & 673.5\tabularnewline 
\hline
\multicolumn{3}{|c||}{} & \multicolumn{5}{|c|}{$\alpha=0.7$, $\alpha'=0.8$}\tabularnewline
\hline 
$k$ & $n$ & $N_\text{sub}$ & $\#$2-level & $n_\text{CS}$ & Time & $\#$1-level & Time\tabularnewline
\hline 
10 & $3.1\times 10^5$ & 125 & 28 & $1.9\times 10^3$ & 8.2 & 58 & 7.7\tabularnewline
15 & $1.5\times 10^6$ & 216 & 39 & $4.2\times 10^3$ & 19.0& 82 & 20.1 \tabularnewline
20 & $6.3\times 10^6$ & 512 & 58 & $7.9\times 10^3$ & 42.4 & 123 & 49.7 \tabularnewline
25 & $1.4\times 10^7$ & 729 & 60 & $1.7\times 10^4$ & 80.6 & 148 & 94.1\tabularnewline 
30 & $3.5\times 10^7$ & 1000 & 80 & $2.6\times 10^4$ & 251.9 & 179 & 328.0\tabularnewline 
\hline
\multicolumn{3}{|c||}{} & \multicolumn{5}{|c|}{$\alpha=0.8$, $\alpha'=0.8$}\tabularnewline
\hline 
$k$ & $n$ & $N_\text{sub}$ & $\#$2-level & $n_\text{CS}$ & Time & $\#$1-level & Time\tabularnewline
\hline 
10 & $3.4\times 10^5$ & 216 & 31 & $1.9\times 10^3$ & 12.6 & 67 & 11.7  \\
20 & $7.1\times 10^6$ & 1000 &70 & $7.9\times 10^3$ &76.9  & 147 & 58.3  \\
30 & $4.1\times 10^7$ & 3375 &109 & $2.6\times 10^4$ & 238.0 & $>$200 & -  \\
40 & $1.3\times 10^8$ & 6859 &193 & $5.1\times 10^4$ & 948.9 & $>$200 & - \\
\hline
\end{tabular}
\caption{$\kappa_\text{prob} = 0$, $\kappa_\text{prec} = k$, $\delta \sim 2h$, $H_\text{sub}  \sim k^{-\alpha}$, $H \sim k^{-\alpha'}$.} 
\label{tab:4}
\end{center}
\end{table}

\bigskip
\noindent \textbf{Acknowledgement} 
This work has been supported in part by the French National Research Agency (ANR), project MEDIMAX, ANR-13-MONU-0012.


\bibliographystyle{plainnat} 
\bibliography{dolean_mini_10}


\end{document}